%
%
%
%
%
%
%
%
%
\scrollmode
\magnification=\magstep1
\parskip=\smallskipamount
\hoffset=1.5cm \hsize=12cm

\def\demo#1:{\par\medskip\noindent\it{#1}. \rm}
\def\ni{\noindent}               
\def\ll{\leftline}
\def\cl{\centerline}

\def\begin{\ll{}\vskip 10mm\nopagenumbers}  
\def\pn{\footline={\hss\tenrm\folio\hss}}   
\def\ii#1{\itemitem{#1}}

%
%
\outer\def\beginsection#1\par{\bigskip
  \message{#1}\leftline{\bf\&#1}
  \nobreak\smallskip\vskip-\parskip\noindent}

%
%
\outer\def\proclaim#1:#2\par{\medbreak\vskip-\parskip
    \noindent{\bf#1.\enspace}{\sl#2}
  \ifdim\lastskip<\medskipamount \removelastskip\penalty55\medskip\fi}

\def\endpr{\hfill $\spadesuit$ \medskip}

%
%
%
%


%
%

\def\R{{\rm I\kern-0.2em R\kern0.2em \kern-0.2em}}
\def\N{{\rm I\kern-0.2em N\kern0.2em \kern-0.2em}}
\def\P{{\rm I\kern-0.2em P\kern0.2em \kern-0.2em}}
\def\B{{\rm I\kern-0.2em B\kern0.2em \kern-0.2em}}
\def\C{{\rm C\kern-.4em {\vrule height1.4ex width.08em depth-.04ex}\;}}
\def\CP{\C\P}

%
%
%
%

\def\cJ{{\cal J}}

\def\cO{{\cal O}}

\def\cR{{\cal R}}
\def\cS{{\cal S}}

%
%
%
\def\a{\alpha}
\def\b{\beta}

\def\d{\delta}
\def\e{\epsilon}

\def\s{\sigma}

\def\c{\chi}

\def\G{\Gamma}

\def\L{\Lambda}

%
%
%
%
\def\bar{\overline}              
\def\bs{\backslash}              
\def\di{\partial}                

%
%

\def\c*{{\C}^*}

\def\wt{\widetilde}

%
%
\def\dim{{\rm dim}}                    
\def\holo{holomorphic}                   
\def\nbd{neighborhood}                   
\def\spsh{strongly\ plurisubharmonic}
\def\ss{\subset\!\subset}                

\def\iff{if and only if}

\def\hvb{holomorphic vector bundle}

\def\Ell{{\rm Ell}}


\begin
\cl{\bf  THE OKA PRINCIPLE FOR SECTIONS OF}
\cl{\bf  SUBELLIPTIC SUBMERSIONS}
\bigskip
\cl{\bf Franc Forstneri\v c}
\bigskip\medskip\rm

%
%
%
%
\beginsection 0. Introduction.

In this paper we prove the Oka principle for maps $X\to Y$ from 
Stein manifolds $X$ to complex manifolds $Y$ which admit a finite dominating 
collection of sprays (such manifolds are called subelliptic), 
as well as the corresponding result for sections of submersions. 
Our main result, Theorem 1.1, extends the results of Oka [O],
Grauert [Gr1, Gr2] and Gromov [G]. We also prove a result on removing 
intersections of \holo\ maps from Stein manifolds with closed complex 
subvarieties with subelliptic complements (section 6).

We begin with a brief survey. In 1939 K.\ Oka [O] proved that a second Cousin problem 
on a domain of holomorphy is solvable if it is solvable by continuous functions.
Oka's result has the following equivalent formulation: 
{\it If $h\colon Z\to X$ is a principal \holo\ fiber bundle with fiber $\C^*=\C\bs \{0\}$
over a domain of holomorphy (or a Stein manifold) then every continuous section 
of $h$ is homotopic to a \holo\ section}. In a seminal work of 1957 H.\ Grauert 
[Gr1, Gr2] proved Oka's theorem with $\C^*$ replaced by any complex Lie group 
or complex homogeneous space, with the stronger conclusion that the inclusion 
${\rm Holo}(X;Z)\hookrightarrow {\rm Cont}(X;Z)$ of the space of 
\holo\ sections into the space of continuous sections is a {\it weak homotopy 
equivalence} (it induces isomorphisms of all homotopy groups of the 
two spaces). This is known as the (parametric) {\bf Oka-Grauert principle}. 
For related results and extensions see [C], [FR] and [HL].

In 1989 M.\ Gromov [G] introduced the concept of a {\it dominating spray} and 
outlined a proof of the {\it parametric Oka principle for sections of  
\holo\ submersions $h\colon Z\to X$ onto a Stein base $X$ which admit 
fiber dominating sprays over small open subsets of $X$} (complete proofs
can be found in [FP1] for fiber bundles and in [FP2, FP3] for submersions). 

In this paper we introduce a more flexible and apparently weaker
condition which also implies the parametric Oka principle. 
Recall that a {\bf spray} on a complex manifold $Y$ is a \holo\ 
map $s\colon E\to Y$ from the total space of a \hvb\ $p\colon E\to Y$ such that 
$s(0_y)=y$ for all $y\in Y$. $Y$ is called {\bf subelliptic} if it 
admits finitely many sprays $s_j\colon E_j\to Y$ such that for any $y\in Y$ the 
vector subspaces $(ds_j)_{0_y}(E_{j,y}) \subset T_y Y$ together span $T_y Y$ 
(Definition 2).  We prove {\it the parametric Oka principle for maps from any 
Stein manifold to any subelliptic manifold, as well as for sections
of subelliptic submersions over a Stein base} (Theorem 1.1).
This extends Gromov's theorem [G, 4.5 Main Theorem] which assumes 
the existence of a {\bf dominating spray} $s\colon E\to Y$ satisfying 
$(ds)_{0_y}(E_y)=T_yY$ for all $y\in Y$ (such manifold $Y$ is called 
{\bf elliptic}). There is no immediate way of creating dominating sprays from 
dominating families of sprays unless the bundles $E_j$ are trivial.

Subellipticity is easier to verify than ellipticity 
and consequently it enables us to extend the Oka principle to 
a wider class of target manifolds. For instance, if $A$ is a closed 
complex ($=$algebraic) subvariety of complex codimension at least two in a 
complex projective space $\CP^n$ (or in a complex Grassmanian) then 
its complement is subelliptic and hence the Oka principle holds for 
maps from Stein manifolds to $\CP^n\bs A$ (Proposition 1.2).
We don't know whether this complement is elliptic in general. 
(By removing a hyperplane we obtain $\C^n\bs A$ which is elliptic [G, FP1].) 
On projective algebraic manifolds subellipticity can be 
localized:  {\it If each point $y\in Y$ admits a Zariski open 
\nbd\ which is algebraically subelliptic then $Y$ is subelliptic} 
(Proposition 1.3). No such result is known about ellipticity. 

Subellipticity is equivalent to the existence of a {\it dominating 
composed spray} (Lemma 2.4). Even though Gromov discussed composed sprays 
in [G] (see in particular the sections 1.3, 1.4.F.\ and 2.9.A.), this condition 
has not been formulated before. On a Stein manifold $Y$ subellipticity is 
equivalent to ellipticity (Lemma 2.2) and both conditions are implied by 
the validity of the Oka principle for maps $X\to Y$ from Stein manifolds 
$X$ with second order interpolation along closed complex submanifolds 
$X_0\subset X$ (see [G, 3.2.A] or [FP3, Proposition 1.2]). 
It is not clear whether the validity of the Oka principle 
implies subellipticity (or ellipticity) for all complex manifolds.  

\pn

\beginsection 1. The results.

Let $h\colon Z\to X$ be a \holo\ submersion onto $X$.
Given a subset $U\subset X$ we write $Z|_U=h^{-1}(U)$. 
For $z\in Z$ we denote by $VT_z Z$ the kernel of $dh_z$
(which equals the tangent space to the fiber $h^{-1}(h(z))$
at $z$) and call it the {\it vertical tangent space} of $Z$ at $z$.
If $p\colon E\to Z$ is a \hvb\ we denote by $0_z\in E$ the 
base point in the fiber $E_z=p^{-1}(z)$. At each point 
$z\in Z$ we have a natural splitting $T_{0_z}E = T_z Z\oplus E_z$.

%
%
%
%
\proclaim Definition 1:  {\rm [G, sec.\ 1.1.B]}
A {\bf spray associated to a \holo\ submersion} $h\colon Z\to X$ 
(an $h$-spray) is a triple $(E,p,s)$, where $p\colon E\to Z$ is 
a \hvb\ and $s\colon E\to Z$ is a \holo\ map such that 
for each $z\in Z$ we have $s(0_z)=z$ and $s(E_z) \subset Z_{h(z)}$.
The spray $s$ is {\bf dominating} at the point $z\in Z$ if the derivative 
$ds \colon T_{0_z} E \to T_z Z$ maps $E_z$ surjectively onto 
$VT_z Z=\ker dh_z$. A {\bf spray on a complex manifold} $Y$ is a spray 
associated to the trivial submersion $Y\to point$. 

%
%
\proclaim Definition 2: A \holo\ submersion $h\colon Z\to X$
is called {\bf subelliptic} if each point in $X$ has an 
open \nbd\ $U\subset X$ such that $h\colon Z|_U \to U$ admits finitely 
many $h$-sprays $(E_j,p_j,s_j)$ for $j=1,\ldots,k$ satisfying  
$$ 
	(ds_1)_{0_z}(E_{1,z}) + (ds_2)_{0_z}(E_{2,z})\cdots 
                     + (ds_k)_{0_z}(E_{k,z})= VT_z Z   \eqno(1.1)
$$
for each $z\in Z|_U$. A collection of sprays satisfying (1.1) is said 
to be {\bf dominating} at $z$. A submersion $h$ is {\bf elliptic} if the above 
holds with $k=1$, i.e., if any point $x\in X$ has a \nbd\ $U\subset X$ 
such that $h\colon Z|_U\to U$ admits a dominating spray.
A complex manifold $Y$ is elliptic (resp.\ subelliptic) 
if the trivial submersion $Y\to point$ is such.

Thus every elliptic submersion is also subelliptic.
Examples of elliptic manifolds and submersions may be found in [G] 
(see especially sections 0.5.B and 3.4.F) and in [FP1].
The exponential map $\exp\colon {\bf g}\to G$ on any 
complex Lie group $G$ gives a dominating spray 
$s\colon E=G\times {\bf g}\to G$, $s(g,t)=\exp(t)g$.

\proclaim 1.1 Theorem: 
If $h\colon Z\to X$ is a subelliptic submersion 
onto a Stein manifold $X$ then the inclusion 
$\iota_h \colon \G_{\rm holo}(X;Z)\hookrightarrow \G_{\rm cont} (X;Z)$ 
of the space of holomorphic sections of $h$ into the space of 
continuous sections is a weak homotopy equivalence. (Both spaces 
are endowed with the topology of uniform convergence on compacts.)

Theorem 1.1 is the main result of this paper. It implies in particular 
that maps $X\to Y$ from any Stein manifold $X$ to any subelliptic manifold $Y$
satisfy the parametric Oka principle (since maps $X\to Y$ correspond to sections 
of the projection $X\times Y\to X$). The proof of Theorem 1.1 will show in addition
that {\it sections of any subelliptic submersion $h\colon Z\to X$ onto a Stein base 
$X$ satisfy the conclusion of Theorem 1.4 in} [FP3] (which is equivalent to the 
$\Ell_\infty$ property introduced by Gromov [G, sec.\ 3.1.]). 
This includes {\it uniform approximation of \holo\ sections on compact 
holomorphically convex subsets of $X$ and interpolation of \holo\ sections 
on any closed complex subvariety $X_0\subset X$}. The extension in [F2] 
to multi-valued sections of ramified mappings $h\colon Z\to X$ onto a 
Stein space $X$ also holds when $h$ is a subelliptic submersion 
outside its ramification locus.

We now give examples of subelliptic manifolds (for proofs see section 5).

\proclaim 1.2 Proposition: (a) If $Y$ is a complex Grassman manifold
and $A \subset Y$ is a closed complex($=$algebraic) subvariety of codimension 
at least two then $Y\bs A$ is subelliptic. This holds in particular 
when $Y$ is a complex projective space $\CP^n$. 
\item{(b)} Let $h \colon Z\to X$ be a \holo\ fiber bundle whose fiber 
is $\CP^n$ or a complex Grassmanian. If $A\subset Z$ is a closed complex 
subvariety whose fiber $A_x=A\cap Z_x$ has codimension at least two in $Z_x$
for any $x\in X$ then the restricted submersion $h\colon Z\bs A \to X$ 
is subelliptic.

The subvarieties in Proposition 1.2 have codimension at least two.
The Oka principle fails in general for maps into complements of  
complex hypersurfaces or non-algebraic subvarieties of any dimension
(see the examples in [FP3]).

Recall that a {\bf  projective algebraic manifold} is a closed 
complex submanifold of a complex projective space and a 
{\bf quasi-projective manifold} is a Zariski open set in a projective 
manifold. We may speak of algebraic vector bundles and {\bf algebraic sprays} 
on such manifolds, and algebraic subellipticity can be localized as follows 
(compare with Lemma 3.5.B.\ and 3.5.C.\ in [G]).

%
%
\proclaim 1.3 Proposition: If $Y$ is a quasi-projective algebraic manifold
such that each point $y\in Y$ has a Zariski open \nbd\ $U\subset Y$ 
and algebraic sprays $s_j\colon E_j\to Y$ ($j=1,2,\ldots,k$), defined
on algebraic vector bundles $p_j\colon E_j\to U$ and satisfying  
$$ 	 
	(ds_1)_{0_y}(E_{1,y}) + (ds_2)_{0_y}(E_{2,y})\cdots 
                     + (ds_k)_{0_y}(E_{k,y})= T_y Y,    
$$ 
then $Y$ is subelliptic.

No such localization result is known for ellipticity. Note that the condition 
in Proposition 1.3 is equivalent to (1.1) when $h$ is the trivial submersion 
$Y\to point$.  The ranges of the sprays $s_j$ need not be contained in $U$.

%
%
\proclaim 1.4 Proposition: If $Y$ is a (quasi-) projective algebraic
manifold with an algebraic spray $s\colon E\to Y$ which is a submersion 
of $E$ onto $Y$ then the complement $Y\bs A$ of every algebraic 
subvariety of codimension at least two is subelliptic.

\proclaim 1.5 Corollary: If $G$ is a complex algebraic Lie group 
whose exponential map is algebraic then the complement $G\bs A$ 
of any algebraic subvariety $A\subset G$ of codimension at least 
two is subelliptic. This holds in particular if $G$ is nilpotent and 
simply connected.

\demo Proof of Corollary 1.5: 
The exponential map $\exp\colon {\bf g}\to G$ is locally biholomorphic and 
hence the spray $s\colon E=G\times {\bf g}\to G$, $s(g,t)=\exp(t)g$, is a 
submersion of $E$ onto $G$. If $G$ and $\exp$ are algebraic 
then $s$ is algebraic and the result follows from Proposition 1.4.
This is the case for simply connected nilpotent Lie groups.
\endpr

We don't know whether all subelliptic manifolds of the form $Y\bs A$ 
considered above (where $A$ is a subvariety of $Y$ containing no hypersurfaces) 
are actually elliptic, although we believe they are not.

%
%

\proclaim 1.6 Proposition: Let $\pi\colon \wt Y\to Y$ be an
unramified \holo\ covering map. If $Y$ is subelliptic (resp.\ 
elliptic) then so is $\wt Y$.

A result of this kind is mentioned in [G, 3.5.B''] 
(see $(**)$ on p.\ 883 of [G]). It is not clear whether the converse 
is true as well, i.e., {\it does (sub-) ellipticity of $\wt Y$ imply 
the same property for $Y$ ?} A good test case may be complex
tori $T=\C^n/\Gamma$ where $\Gamma\subset \C^n$ is a lattice of 
real rank $2n$. Denote by $\pi\colon \C^n\to T$ the quotient map
(which is a universal covering of $T$). The spray 
$s\colon \C^n\times\C^n\to\C^n$, $s(z,t)=z+t$, is 
$\Gamma$-equivariant and hence it passes down to a spray on  
$T$. Removing the point $p_0=\pi(0)\in T$ we obtain 
a covering map $\pi \colon \C^n\bs \Gamma \to T\bs \{p_0\}$.
It is easily seen that for $n\ge 2$ the lattice $\Gamma$ is a tame 
discrete subset of $\C^n$ in the sense of Rosay and Rudin [RR]
and hence $\C^n\bs \Gamma$ admits a dominating spray according
to Lemma 7.1 in [FP2]. However, these sprays don't pass to
sprays on $T\bs \{p_0\}$ and it is not clear whether the latter
manifold is subelliptic.

%
%
\ni\bf Outline of the paper. \rm Section 2 contains some basic results 
and constructions with sprays. We show that the domination of a family 
of sprays is equivalent to the domination of the associated {\it composed spray}. 
In section 3 we prove a homotopy version of the Oka-Weil approximation theorem 
for sections of submersions onto a Stein base which admit a finite dominating 
collection of sprays (Theorem 3.1). In section 4 we prove Theorem 1.1. In section 5 
we prove Propositions 1.2, 1.3, 1.5 and 1.6. In section 6 we use the methods 
developed in section 4 (and in [FP2], [F1]) to prove a result on removing 
intersections of \holo\ maps $X\to Y$ from Stein manifolds $X$ with any closed complex 
subvariety $A\subset Y$ whose complement $Y\bs A$ is subelliptic. Theorem 6.1 
contains as special cases the well known theorem of Forster and Ramspott 
on complete intersections [FR] as well as Theorem 1.3 from [F1].
The appendix contains some remarks on Gromov's paper [G].

%
%
%
%
\beginsection 2. Subellipticity and composed sprays.

In this section we first collect some basic results on subelliptic 
manifolds and submersions. The constructions in Lemmas 2.1 and 2.3 are 
due to Gromov [G]. Recall that a finite collection of sprays 
$s_j\colon E_j\to Y$ ($j=1,2,\ldots,k$) is {\it dominating} if the 
subspaces $(ds_j)_{0_y}(E_{j,y}) \subset T_y Y$ together span $T_y Y$ 
for each $y\in Y$ (Definition 2).

\proclaim 2.1 Lemma: If $s_j\colon E_j\to Y$ ($1\le j\le k$)
is a dominating collection of sprays on $Y$, defined on trivial 
bundles $E_j\simeq Y\times \C^{m_j}$, then $Y$ admits
a dominating spray. The analogous result holds for $h$-sprays.

\demo Proof:  We may assume that $s_j$ is defined on $Y\times \C^{m_j}$
for each $j=1,2,\ldots,k$. We define sprays 
$s^{(j)} \colon Y\times\C^{m_1+\ldots+m_j} \to Y$ inductively 
by $s^{(1)}=s_1$ and 
$$ 
	s^{(j)}(y,e_1,\ldots,e_j) =
	s_j\bigl( s^{(j-1)}(y,e_1,\ldots,e_{j-1}), e_j\bigr),\quad 2\le j\le k.   
$$
Clearly we have 
$$ (ds^{(k)})_{0_y} (E^{(k)}_y) = 
   (ds_1)_{0_y}(E_{1,y}) + (ds_2)_{0_y}(E_{2,y})\cdots 
                     + (ds_k)_{0_y}(E_{k,y}).
$$
Hence $(s_1,\ldots,s_k)$ is a dominating collection of sprays
on $Y$ \iff\ $s^{(k)}$ is a dominating spray. We call $s^{(k)}$
the {\it direct sum} of the sprays $s_j$ and write
$s^{(k)} = s_1\oplus\cdots\oplus s_k$. (Observe that this
construction is possible only for sprays defined on trivial bundles.)
\endpr

%
%
\proclaim 2.2 Lemma: Any subelliptic Stein manifold is elliptic. 
If $Z,X$ are Stein manifolds then any subelliptic submersion 
$h\colon Z\to X$ is elliptic.  

\demo Proof: \rm Assume that $Y$ is Stein and $s_j\colon E_j\to Y$ ($j=1,\ldots,k$) 
is a dominating collection of sprays defined on vector bundles 
$p_j\colon E_j\to Y$. By Cartan's Theorem A [GR]  any \hvb\ over $Y$ is generated 
by finitely many global \holo\ sections. This gives for each $j$ a surjective 
vector bundle map $g_j\colon Y\times \C^{m_j}\to  E_j$ for some 
large $m_j\in \N$. Then $\wt s_j= s_j\circ g_j \colon Y\times \C^{m_j} \to Y$ is a 
spray whose vertical derivative at the zero section has the same range as that 
of $s_j$. It follows that $\wt s_1\oplus \cdots \oplus \wt s_k$ (defined in Lemma 2.1)
is a dominating spray on $Y$. A similar proof holds for submersions.
\endpr

For general sprays one can use the following construction.

\proclaim Definition 3: {\rm (Gromov [G, sec.\ 1.3])} 
Let $(E_1,p_1,s_1)$ and $(E_2,p_2,s_2)$  be $h$-sprays associated to a 
\holo\ submersion $h\colon Z\to X$. The {\bf composed spray} 
$(E_1\# E_2, p_1\# p_2, s_1\# s_2)=(E^*,p^*,s^*)$ 
is defined by
$$ \eqalign{ & E^* = 
     \{(e_1,e_2)\in E_1\times E_2 \colon\ \ s_1(e_1)=p_2(e_2) \}, \cr
   & \ \ p^*(e_1,e_2) = p_1(e_1), \quad s^*(e_1,e_2) = s_2(e_2). \cr}
$$

This operation extends to any finite collection of sprays
and it includes iterations of sprays (and of composed sprays). 
The $k$-th iteration $(E^{(k)},p^{(k)},s^{(k)})$ of $(E,p,s)$ is 
$$ 
   \eqalign{ & E^{(k)} = \{e=(e_1,e_2,\ldots,e_k) \colon e_j\in E\
            \ {\rm for}\ j=1,2,\ldots,k, \cr
            &\qquad \qquad\qquad
              s(e_j)=p(e_{j+1})\ {\rm for}\ j=1,2,\ldots,k-1\},  \cr
            & p^{(k)}(e) = p(e_1),\quad s^{(k)}(e) = s(e_k).  \cr}
$$
Composed sprays are not sprays in the sense of Definition 1 since $E_1\# E_2$ 
does not have a natural structure of a \hvb\ over $Z$. Observe that $E_1\# E_2$ 
is the pull-back $s_1^*(E_2)$ of the vector bundle $p_2\colon E_2\to Z$ 
by the first spray map $s_1\colon E_1\to Z$ (and hence is a \hvb\ over the 
total space of the bundle $E_1\to Z$). However, the composed bundle has a well 
defined zero section and a partial linear structure on  fibers.

\proclaim 2.3 Lemma: If $(E_j,p_j,s_j)$ for $j=1,\ldots,k$ are $h$-sprays 
on $Z$ then the restriction of the composed bundle $E_1\#\cdots\#E_k\to Z$ 
to any Stein subset $\Omega\subset Z$ is isomorphic to the direct sum
bundle $E_1\oplus E_2\oplus\cdots \oplus E_k|_\Omega$.
Explicitly, there exists a biholomorphic map 
$\theta\colon \oplus E_j|_{\Omega} \to \# E_j|_{\Omega}$ which maps
fibers to fibers and preserves the zero section. (The set $\Omega$ may be 
either an open Stein subset or a Stein subvariety of $Z$.)

\demo Proof: It suffices to consider the case $k=2$ and apply induction. 
By the construction the composed bundle $E=E_1\# E_2$ is a \hvb\ over $E_1$ 
with projection $p\colon E\to E_1$. The total space $E'_1=E_1|_{\Omega}$ 
of the restricted bundle $p_1\colon E_1|_{\Omega}\to \Omega$ is a Stein manifold.
Let $E'=E|_{E'_1}$. Applying Grauert's theorem [Gr2] to the homotopy 
$$ g_t\colon E'_1\to E'_1, \quad g_t(z,e)=p_1(z,te), 
					\qquad (z,e)\in E'_1,\ t\in [0,1] 
$$ 
we see that the pull-backs $g_t^*(E')\to E'_1$ of $p\colon E'\to E'_1$
are all isomorphic. Since $g_1$ is the identity on $E'_1$ 
and $g_0=p_1$, this gives an isomorphism between $p\colon E'\to E'_1$ 
and $p_1^*(E|_{\Omega})$, where $\Omega$ denotes the zero section of $E'_1$. 
Since $s_1$ is the identity on the zero section  of $E_1$, we have 
$E|_{\Omega}=E_2|_{\Omega}$. Hence the bundle $p\colon E'\to E'_1$ 
is isomorphic to $\pi\colon p_1^*(E_2|_{\Omega}) \to E'_1$ 
(as a bundle over $E'_1$). By linear algebra the composition 
$p_1\circ \pi\colon p_1^*(E_2|_{\Omega})\to \Omega$ is isomorphic 
to $E_1\oplus E_2|_{\Omega}$ (as a \hvb\ over $\Omega$), and this
endows $p_1\circ\pi \colon E'\to\Omega$ with the structure of a \hvb\ 
over $\Omega$ isomorphic to $E_1\oplus E_2|_{\Omega}$. 
\endpr

The notion of being `dominating' (Definition 1) carries over in an obvious 
way to composed sprays by requiring the submersivity of the spray map
in the fiber direction along the zero section of the composed bundle. 
The next lemma  follows immediately from definitions and explains the 
relevance of subellipticity.

%
%
\proclaim 2.4 Lemma: Let $h\colon Z\to X$ be a \holo\ submersion.
A collection of $h$-sprays $s_j\colon E_j\to Z$ $(j=1,2,\ldots,k)$ 
is dominating at $z\in Z$ \iff\ the composed $h$-spray 
$s_1\#\cdots\#s_k\colon E_1\#\cdots\#E_k\to Z$ is dominating at $z$. 

For the sake of completeness we also add the following result.
\proclaim 2.5 Lemma: The Cartesian product of any finite family 
of (sub-) elliptic manifolds is (sub-) elliptic.

\demo Proof: It suffices to prove the result for the product
of two manifolds. Let $Y=Y_1\times Y_2$ and let
$\pi_j\colon Y\to Y_j$ ($j=1,2$) denote the projection $\pi_j(y_1,y_2)=y_j$. 
If $(E_j,p_j,s_j)$ is a spray on $Y_j$ for $j=1,2$, we get a spray $s=s_1\oplus s_2$ 
on the bundle $E=\pi_2^*E_1\oplus \pi_1^*E_2 \to Y$ given by
$$ 
	s(y_1,y_2,e_1,e_2) = \bigl(s_1(y_1,e_1), s_2(y_2,e_2) \bigr).
$$
If $s_1$ is dominating on $Y_1$ and $s_2$ is dominating on $Y_2$ then
$s$ is dominating on $Y$. Similarly, if a family of sprays
$\{s_i\colon i=1,\ldots,i_0\}$ is dominating on $Y_1$ and a family 
of sprays $\{\sigma_k\colon k=1,\ldots,k_0\}$ is dominating on $Y_2$ then
the collections $s_i\oplus \sigma_k$ is dominating on $Y$.
\endpr

%
%
%
\beginsection 3. The Oka-Weil theorem for subelliptic submersions.

%
%
%
\proclaim 3.1 Theorem:
Let $h\colon Z\to X$ be a holomorphic submersion of a complex manifold
$Z$ onto a Stein manifold $X$. Assume that there exist $h$-sprays 
$s_j\colon E_j\to Z$ ($j=1,2,\ldots,k$) which together dominate at each 
point $z\in Z$ (i.e., condition (1.1) holds). Let $d$ be a metric 
on $Z$ inducing the manifold topology. Suppose that $K$ is a compact 
holomorphically convex set in $X$, $U\supset K$ is an open set and 
$f_t\colon U \to Z$ $(t\in [0,1])$ is a homotopy of \holo\ sections such that 
$f_0$ extends to a holomorphic section on $X$. Then for any $\e>0$ there exists a 
homotopy of holomorphic sections $\wt f_t\colon X\to Z$ such that $\wt f_0=f_0$ and
$d\bigl(\wt  f_t(x),f_t(x)\bigr) < \e$ for $x\in K$ and $(t\in [0,1])$.
The analogous result holds for parametrized families of sections.

\demo Proof: For submersions which admit a global dominating $h$-spray
this is Theorem 4.1 in [FP1], and the general parametric case is Theorem 4.2 
in [FP1]. (See also Theorems 2.1 and 2.2 in [FP3] for the Oka-Weil theorem 
with interpolation on a complex subvariety in $X$.) To prove Theorem 3.1 
we replace the dominating $h$-spray $s\colon E\to Z$ in the proof of 
Theorem 4.1 in [FP1] by the doninating composed $h$-spray 
$$ 
     s=s_1\# s_2\#\cdots\# s_k \colon E = E_1\#E_2\cdots\#E_k \to Z 
$$
(Definition 3). The main point in the proof of Theorem 4.1 in 
[FP1, p.\ 135] is that for some sufficiently large $k\in\N$ 
(depending only on $\{f_t\}$) there exists a homotopy of \holo\ sections 
$\xi_t$ ($t\in [0,1]$) of the iterated spray bundle $E^{(k)}|_{f_0(V)} \to Z$, 
restricted to $f_0(V)\subset Z$ for a sufficiently small open set $V\supset K$, 
such that 
$$
   s^{(k)} \bigl(\xi_t(f_0(x))\bigr) = f_t(x),\quad (x\in V,\ t\in [0,1]). \eqno(3.1)
$$ 
The same is true in the present situation which can be seen as follows.
Since $s\colon E|_{f_0(V)}\to Z$ is a submersion on the zero section, 
there are a number $t_1>0$ (depending only on $\{f_t\}$) and a homotopy 
$\xi_t$ of \holo\ sections of the latter bundle satisfying (3.1) 
for $k=1$ and $0\le t\le t_1$. Repeating the argument (and shrinking 
$V \supset K$ if necessary) we obtain a number $t_2>t_1$ 
(depending only on $\{f_t\}$) and a family of sections 
$\{\xi_t \colon t_1\le t\le t_2\}$ of $E|_{f_{t_1}(V)}$ such that 
$\xi_{t_1}$ is the zero section and  
$$ 
  s\bigl( \xi_t(f_{t_1}(x)) \bigr) = f_t(x),\quad (x\in V,\ t\in[t_1,t_2]).
$$
The two homotopies $\xi_t$ together for $0\le t\le t_2$ define a homotopy of 
sections of the second iteration $E^{(2)}|_{f_0(V)}$ such that (3.1) holds for 
$k=2$ and $t\in [0,t_2]$. In finitely many steps (whose number depends only on 
$\{f_t\}$) we obtain a family $\xi_t$ satisfying (3.1).

By Lemma 2.3 $E^{(k)} \to Z$ admits the structure of a \hvb\ over any 
Stein subset of $Z$. Hence the usual Oka-Weil theorem holds for sections 
of $E^{(k)}|_{f_0(X)} \to f_0(X)$. This gives a homotopy of \holo\ sections 
$\wt\xi_t$ of $E^{(k)}|_{f_0(X)}$ for $t\in[0,1]$ such that $\wt\xi_t$ 
approximates $\xi_t$ uniformly on $f_0(K)$ for each $t\in[0,1]$ and $\wt\xi_0$ 
is the zero section. The homotopy 
$$
	\wt f_t(x)=s^{(k)}(\wt\xi_t(f_0(x))) \qquad (x\in X,\ t\in [0,1])
$$ 
of sections of $h\colon Z\to X$ then satisfies Theorem 3.1.
Similarly one proves the parametric version of Theorem 3.1 
(see Theorem 4.2 in [FP1]).

\beginsection 4. Proof of Theorem 1.1.

We shall follow the proof of Theorem 1.5 in [FP2] (or Theorem 1.4 in [FP3])
and explain the necessary modifications. The basic problem is to deform a 
continuous section of $h\colon Z\to X$ to a \holo\ section. The construction 
proceeds through a sequence of modifications in which we obtain \holo\ sections 
over a family of open subsets which increase to $X$. The proof has two 
essential ingredients:

\item{(1)} Solution of the {\it modification problem} explained below, and  

\item{(2)} An inductive construction of a sequence of {\it holomorpic complexes} 
(in the terminology of [FP2, FP3]) which converges uniformly on compacts in $X$
to a global \holo\ section. 

\ni Part (2) (globalization) uses the solution of part 1 and is explained in 
section 5 of [FP2] or (with slightly less details) in [FP3]). This part does 
not require any changes whatsoever. The rest of this section is devoted to
part 1. 

An ordered pair $(A,B)$ of compact subsets of $X$ is said to be a 
{\bf Cartan pair} in $X$ (Definition 4.1 in [FP2]) if 
\item{(i)} $A$, $B$, and $A\cup B$ have a basis of Stein \nbd s,
\item{(ii)} $\bar{A\bs B} \cap \bar{B\bs A} =\emptyset$, and 
\item{(iii)}  the set $C=A\cap B$ is Runge in $B$.
($C$ may be empty.)

\ni\bf The modification problem. \rm Let $(A,B)$ be a Cartan pair in $X$
and let $a\colon \wt A\to Z$, $b\colon \wt B\to Z$ be \holo\ sections of 
$h\colon Z\to X$ in open \nbd s $\wt A\supset A$ resp.\ $\wt B\supset B$. 
Suppose furthermore that $b_t$ ($t\in[0,1]$) is a family of \holo\ sections 
over $\wt C=\wt A\cap \wt B$ such that $b_0=b|_{\wt C}$ and $b_1=a|_{\wt C}$.
Assume that $h\colon Z|_{\wt B}\to \wt B$ admits a finite dominating family
of $h$-sprays. The goal is to construct a \holo\ section $\wt a$ over a \nbd\ of 
$A\cup B$ which is uniformly close to $a$ on $A$ and is obtained 
from the initial sections $a,b$ by homotopies over small \nbd s 
of $A$ resp.\ $B$. (The analogous modification problem must be 
considered for parametrized families of sections as in [FP1]. However, 
it will suffice to explain the proof of the basic non-parametric case 
since the general case follows the same pattern as in [FP1].) 

The modification problem will be solved as in [FP1] by performing 
the following steps.

\ni{\bf Step 1 (approximation):} Find a family of \holo\ sections
$\wt b_t$ ($t\in [0,1]$) over a \nbd\ of $B$ such that 
$\wt b_0=b$ and $\wt b_t$ approximates $b_t$ over an open
\nbd\ of $C$ for each $t\in [0,1]$.

\ni{\bf Step 2 (gluing):} Replace $b$ by $\wt b_1$ from Step 1
and `glue' the pair of sections $a,b$ (which are uniformly close over 
a \nbd\ of $C$) into $\wt a$.

Step 1 has been accomplished by Theorem 3.1 in the present paper 
(which replaces Theorem 4.1 in [FP1]). Step 2 is accomplished by 
the following result which replaces Theorems 5.1 and 5.5 in [FP1]. 

%
%
%
%
\proclaim 4.1 Theorem:
Let $h\colon Z\to X$ be a holomorphic submersion onto a Stein
manifold $X$ and let $d$ be a metric on $Z$ compatible with the 
manifold topology. Let $(A,B)$ be a Cartan pair in $X$. 
Assume that $\wt B\supset B$ is an open set in $X$ such that 
$h\colon Z|_{\wt B}\to \wt B$ admits a finite dominating 
family of $h$-sprays (Definition 2). Let $a\colon \wt A\to Z$ be a 
\holo\ section in an open set $\wt A\supset A$. Then for each $\e>0$ 
there is a $\d>0$ satisfying the following. If $b\colon \wt B\to Z$ 
is a \holo\ section satisfying $d\bigl( a(x), b(x)\bigr) <\d$ for 
$x\in \wt A\cap \wt B$, there exists a homotopy $a_t$ (resp.\ $b_t$) 
of holomorphic sections over a \nbd\ $A'$ of $A$ (resp.\ over a \nbd\ $B'$ 
of $B$) such that $a_0=a|_{A'}$, $b_0=b|_{B'}$, 
$a_1|_{C'}=b_1|_{C'}$ for $C'=A'\cap B'$, and
$$ \eqalign{ d\bigl( a_t(x), a(x) \bigr) &< \e \quad (x\in A',\ t\in [0,1]); \cr
    d\bigl( b_t(x), b(x) \bigr) &< \e \quad (x\in C',\ t\in [0,1]). \cr}
$$
The analogous result holds for parametrized families of
sections, i.e., the conclusion of Theorem 5.5 in [FP1] holds
in the present context.

\demo Proof of Theorem 4.1: 
We shall only prove the basic non-parametric case (for the parametric case we refer 
to Theorem 5.5 in [FP1]). Under the stronger assumption that $h$ admits a dominating 
spray over $\wt B$ this is Theorem 5.1 in [FP1] which was reduced to Proposition 5.2 
in [FP1] (the model case) by using Lemma 5.4 in [FP1]. Unfortunately the proof of 
this lemma does not hold under the current weaker assumption. The following result 
is a suitable replacement. We denote by $\B^n(\e)$ the open ball with radius 
$\e$ in $\C^n$ with center at the origin.  

\proclaim 4.2 Lemma: Let $U$ be an open Stein subset of $Z$
and $s_j\colon U\times \B^n(\e)\to Z$ $(j=1,2)$ \holo\ submersions
such that $s_j(z,0)=z$ and $h(s_j(z,t))=h(z)$ for $z\in U$, $t\in \B^n(\e)$.
Let  $M_j=\{(z,t)\in U\times\C^n\colon (ds_j)_{0_z}t=0\}$.
If the bundles $M_1\to U$ and $M_2\to U$ are isomorphic then 
for any relatively compact subset $V\ss U$ there exist numbers 
$\d>0,\eta>0$, with $0<\eta<\e$, satisfying the following. For any pair 
of points $z,w\in V$ with $h(z)=h(w)$ and $d(z,w)<\d$ there exists 
an injective map $\phi(z,w,\cdotp)\colon \B^n(\eta)\to \B^n(\e)$ 
which is \holo\ in all variables and satisfies 
$$ 
	s_2(w,\phi(z,w,t))=s_1(z,t),\quad \phi(z,z,0)=0.
$$  
If $s_2$ is uniformly close to $s_1$ then we may choose $\phi$ 
such that $\phi(z,z,t) \approx t$.

\demo Proof: Since $U$ is Stein, there is a \holo\ splitting 
$U\times\C^n=M_j\oplus N_j$ for some holomorphic vector subbundle 
$N_j \subset U\times \C^n$. The differential of $s_j$ at the zero section 
carries $N_j$ isomorphically onto $VT(Z)|_U$ and hence 
$N_1\simeq N_2\simeq VT(Z)|_U$. Since $M_1\simeq M_2$ by assumption, 
there exists an automorphism $\theta$ of the trivial bundle $U\times\C^n$ 
with $\theta(M_1)=M_2$ and $\theta(N_1)=N_2$. Set 
$\wt s_2=s_2\circ\theta\colon U\times\C^n\to Z$. The kernel
of $d\wt s_2$ at the zero section equals $M_1$. If we can find
a map $\wt \phi$ satisfying the conclusion of the lemma for  $s_1,\wt s_2$ 
then the map $\phi$ defined by $\theta(w,\wt \phi(z,w,t))=(w,\phi(z,w,t))$
satisfies it for $s_1,s_2$.

Hence we may assume that $M_1=M_2=M$ and $N_1=N_2=N$.
We split the fiber vectors $t=(t',t'') \in M_z\oplus N_z=\C^n$ 
accordingly (the splitting depends on the base point $z\in U$). 
The inverse function theorem shows that for each $z\in U$ the restriction 
$s_1(z,0',\cdotp) \colon N_z \to Z_{h(z)}$ maps a \nbd\ of $0''$ 
in $N_z$ biholomorphically onto a \nbd\ of $z$ in the fiber $Z_{h(z)}$. 
The same is true for the map
$$  t''\in N_z \to  s_1(z,t',t'')\in Z_{h(z)}                   \eqno(4.1) $$
for all sufficiently small $t'\in M_z$. If $w\in Z_{h(z)} \cap U$ is 
chosen sufficiently close to $z$ and if $t'\in M_z$ is sufficiently small 
then by the same argument the map 
$$  t''\in N_z \to   s_2(w,t',\cdotp) \in Z_{h(z)}              \eqno(4.2) $$ 
maps a \nbd\ of $0''$ in $N_z$ biholomorphically onto a \nbd\ of $w$ in 
the fiber $Z_{h(w)}=Z_{h(z)}$. If $w$ is chosen sufficiently close to $z$,
the image of the latter \nbd\ also contains the point $z$ and we let  
$\phi'(z,w,t',\cdotp) \colon N_z \to N_z$ be the composition of (4.1) 
with the (unique) local inverse of (4.2). The map 
$$ 
    (t',t'')\in M_z\oplus N_z\to \phi(z,w,t',t'')= 
    \bigl(t', \phi'(z,w,t',t'') \bigr), 
$$
which is defined for all sufficiently small $t=(t',t'')\in M_z\oplus N_z$,
satisfies Lemma 4.2. 
\endpr

We continue with the proof of Theorem 1.1. Write $C=A\cap B$. 
By Lemma 5.3 in [FP1] there exists a Stein open set $V\subset Z$ containing 
$a(A)$ and a \holo\ submersion (a local spray) $s\colon V\times \B^n(\eta)\to Z$ 
for some $\eta>0$ and $n\in\N$ such that $s(z,0)=z$ and $h(s(z,t))=h(z)$. 
(It is important that $s$ is a submersions over a \nbd\ of $a(C)$.)

By assumption there exists a dominating family of $h$-sprays  
$(E_j,p_j,\s_j)$ ($j=1,\ldots, k$) over $Z|_{\wt B}=h^{-1}(\wt B)$.
Set $E=E_1\oplus\cdots\oplus E_k$ and $\wt E=E_1\#\cdots\# E_k$ (section 2).
Let $\wt\s\colon \wt E\to Z|_{\wt B}$ denote the composed spray which is 
dominating by Lemma 2.4. Choose a Stein open set $U\subset Z$ with 
$a(C) \subset U\subset V\cap (Z|_{\wt B})$. Lemma 2.3 gives a fiber preserving 
biholomorphic map $\theta\colon E|_U\to \wt E|_U$. The map 
$\wt \s\theta \colon E|_U \to Z$ is then a dominating $h$-spray defined 
on the vector bundle $E|_U$, with values in $Z$. 

Set $v_j(z)=(ds)_{0_z}(\di/\di t_j)$ for $z\in V$ and $j=1,\ldots, n$, 
where $s$ is a local spray on $V\supset a(A)$ as above and 
$t=(t_1,\ldots,t_n)\in\C^n$. There exist \holo\ sections $\wt v_j$ of 
the bundle $E|_U \to U$ such that $d(\wt \s\theta)_0(\wt v_j)=v_j$ for $j=1,\ldots,n$. 
(Such sections are unique in any \holo\ subbundle $N\subset E|_U$ which is 
complementary to the kernel of the differential of $\wt \s\theta$ at the zero section.) 
Let $\tau\colon U\times\C^n\to E|_U$ be the map 
$\tau(z,t_1,\ldots,t_n)=\sum_{j=1}^n t_j \wt v_j(z) \in E_z$. The 
composition $\s=\wt \s\theta\tau\colon U\times\C^n\to Z$ is then
a dominating $h$-spray on $U$ (with values in $Z$) satisfying 
$(d\s)_{0_z}(\di/\di t_j)=v_j(z)=(ds)_{0_z}(\di/\di t_j)$ for 
$z\in U$ and $j=1,\ldots,n$. By Lemma 4.2 there are $\d,\eta>0$ 
and a \holo\ map $\phi$ satisfying
$$ 
	\s(w,\phi(z,w,t))=s(z,t),\quad \phi(z,z,0)=0
$$ 
for $t\in \B^n(\eta)$ and $z,w\in U$ with $h(z)=h(w)$ and $d(z,w)<\d$.

Let $b\colon \wt B\to Z$ be a \holo\ section satisfying $d(a(x),b(x))<\d$ 
for $x\in \wt A\cap \wt B$. By decreasing $\d$ we may assume $b(C)\subset U$. 
The problem is that the spray $\s$ is only defined over $U$ and 
need not extend holomorphically into any \nbd\ of $b(B)$. To complete
the proof we shall approximate $\s$ uniformly in a \nbd\ of $b(C)$ by a spray 
$\s'$ which is \holo\ in \nbd\ of $b(B)$. We proceed as follows.

Since $C$ is Runge in $B$ and both sets have a basis of Stein \nbd s in $X$, 
there exist open Stein sets $U'\subset W$ in $Z$ such that
$U'$ is Runge in $W$, $b(B)\subset W$, and $b(C)\subset U' \subset U\cap W$.
Fix a compact set $K\subset E|_{U'}$ and let $K'\subset U'$ denote its base projection.
By Lemma 2.3 and the Oka-Weil theorem there is a fiber preserving holomorphic 
map $\theta'\colon E|_W\to \wt E|_W$ which is uniformly close to 
$\theta \colon E|_U\to \wt E|_U$ on $K$. Also there is for each $j=1,\ldots, n$ 
a \holo\ section $v'_j$ of $E|_{W}$ which approximates $\wt v_j$ on $K'\subset U'$.
Define $\tau'\colon W\times\C^n\to E|_W$ by $\tau'(w,t)=\sum_{j=1}^n t_j v'_j(w)$. 
Then the spray $\s'=\wt \s\theta'\tau'\colon W\times \C^n\to Z$ is 
uniformly close to $\s=\wt \s \theta\tau$ on the compact set 
$L=(\tau')^{-1}(K) \subset U'\times\C^n$. For each $U_1\ss U'$ we may 
choose $K$ sufficiently large to insure that $U_1\times \B^n(3) \subset L$.

Applying Lemma 4.2 to the pair of sprays $\s$, $\s'$ we obtain 
(after shrinking $U'\supset b(C))$ a \holo\ map 
$\xi\colon U'\times \B^n(2)\to \C^n$ which is uniformly close to 
$\xi_0(w,t)=t$ on $U_1 \times \B^n(1)$ and which satisfies 
$\s'(w,\xi(w,t))=\s(w,t)$ and $\xi(w,0)=0$ ($w\in U'$, $t\in \B^n(2)$). Then
$$  
	\s'(w,\xi(w,\phi(z,w,t)))= \s(w,\phi(z,w,t))= s(z,t)
$$
for any $t\in\B^n(\eta)$ and any pair of point $z\in U$, $w\in U'$ 
with  $h(z)=h(w)$. Choose open sets $A', B'\subset X$ satisfying 
$A\subset A'\ss \wt A$, $B\subset B'\ss \wt B$, $C'=A'\cap B'$,
$\bar{a(A')}\subset V$, $\bar{b(B')} \subset W$, $\bar{b(C')} \subset U'$. Set 
$$ 
	\eqalign{ s_1(x,t) &= s(a(x),t),   \qquad\  (x\in A',\ t\in \B^n(\eta)); \cr
             s_2(x,t) &= \s'(b(x),t), \qquad (x\in B',\ t\in\C^n); \cr
            \psi(x,t) &= \xi(b(x),\phi(a(x),b(x),t)), 
                         \quad (x\in C',\ t\in \B^n(\eta)). \cr}
$$
Then $s_2(x,\psi(x,t))=s_1(x,t)$ for $x\in C'$ and $t\in \B^n(\eta)$.
If $b$ is uniformly close to $a$ on $\wt C$ then (since $\xi$ is close 
to $\xi_0(w,t)=t$) the map $\psi$ is uniformly close to 
$\psi_0(x,t)=\phi(a(x),a(x),t)$ ($x\in C'$, $t\in \B^n(\eta)$).
Note that $\psi_0(x,0)=0$ for all $x\in C'$. 

We have thus reduced Theorem 4.1 to Proposition 5.2 in [FP1]. 
To complete the proof we shrink the sets $A'\supset A$, $B'\supset B$
and take $\a\colon A'\to \B^n(\eta)$, $\b\colon B'\to \C^n$ to be \holo\ maps 
furnished by that proposition, satisfying
$\psi(x,\a(x))=\b(x)$ for $x\in C'=A'\cap B'$. Then 
$$
	s_2(x,\b(x))=s_2(x,\psi(x,\a(x)))=s_1(x,\a(x)),\quad (x\in C')
$$ 
and hence these expressions define a \holo\ section $\wt a\colon A'\cup B'\to Z$.
Further details can be found in [FP1].
\endpr

\ni\bf A remark on [FP3]. \rm The condition $M_1\simeq M_2$ (i.e., the
kernels of $ds_1$ and $ds_2$ along the zero section are isomorphic) is 
necessary for the existence of a map $\psi$ as above. This was not stated 
explicitly in the proof of Theorem 5.1 in [FP3]. However, the construction 
in [FP1] (see especially Lemma 5.4 in [FP1]) produces a pair of submersions 
$s_1$, $s_2$ for which this condition is satisfied, and one can apply the
proof in [FP3] to such a pair.
\endpr

%
%
%
%
\beginsection 5. Subelliptic manifolds.

In this section we prove Propositions 1.2, 1.3, 1.5 and 1.6.

\demo Proof of Proposition 1.2: 
In part (a) we shall give a detailed calculation only for $Y=\CP^n$ and will
observe that the same proof applies to complex Grassmanians. Similar arguments
apply to part (b) and we omit the details (compare with the proof of Corollary 1.8 in [FP2]). 

Given an algebraic subvariety $A\subset \CP^n$ containing no complex hypersurfaces
we wish to construct a finite dominating family of algebraic sprays on $\CP^n\bs A$.  
We begin by choosing a hyperplane $\Lambda\subset \CP^n$ and homogeneous coordinates
$Z=[Z_0\colon Z_1\colon \ldots\colon Z_n]$ on $\CP^n$ such that $\Lambda=\{Z_0=0\}$.
Set $U_j=\{Z\in\CP^n\colon Z_j\ne 0\} \simeq\C^n$ for $j=0,1\ldots,n$ (hence
$\CP^n=U_0\cup \Lambda$).  Let $L\to \CP^n$ denote the \holo\ line bundle 
$L=[\Lambda]^{-1}$ where $[\Lambda]$ is the line bundle determined by 
the divisor of $\Lambda$. (The usual notation is $L=\cO_{\CP^n}(-1)$, see [GH].)
$L$ admits \holo\ trivializations $\phi_j\colon L|_{U_j}\to U_j\times \C$ 
with transition maps
$$ 
	\phi_{ik}(Z,t)=\phi_i\circ \phi^{-1}_{k}(Z,t)= (Z,tZ_i/Z_k),\quad  
	(Z\in U_{ik}=U_i\cap U_k,\ t\in \C).
$$
Choose a vector $v\in\C^n$ such that the orthogonal projection 
$\pi \colon U_0=\C^n\to\C^{n-1}$ with kernel $\C v$ is proper when restricted 
to $A\cap U_0$. (This is the case precisely when the complex line $\C v$
does not intersect $A$ at any point of $\Lambda$. Since $A$ has codimension
at least two,  it does not contain $\L$ and hence this holds for almost
every $v$.) Then $A'=\pi(A\cap U_0)\subset\C^{n-1}$ is a proper algebraic
subvariety of $\C^{n-1}$. Let $p$ be any nonzero \holo\ polynomial on
$\C^{n-1}$ which vanishes on $A'$. Then the map 
$U_0\times \C\to U_0$ given by
$$ 
	s(z,t)=  z+tp(\pi z)v =z+tf(z) \qquad (z\in U_0,\ t\in \C) 
$$
is a spray on $U_0$ with ${\di\over \di t}s(z,0)=p(\pi z)v=f(z)$. Although
$s$ does not extend to a spray $\CP^n\times\C\to\CP^n$ because of
singularities in $\Lambda=\{Z_0=0\}$, it induces a spray 
$\wt s\colon L^{\otimes m} \to \CP^n$ where $m$ is the degree of the
polynomial $p$ and $E=L^{\otimes m}$ denotes the $m$-th tensor 
power of $L$. $E$ admits trivializations 
$\theta_i\colon E|_{U_i}\to U_i\times\C$ ($i=0,1,\ldots,n$)
with transition maps
$$  
    \theta_{ik}(Z,t) = \bigl( Z,t(Z_i/Z_k)^m \bigr),\qquad  
    (Z\in U_{ik},\ t\in \C).                           
$$
Set $\wt s = s\, \wt\theta_0 \colon E|_{U_0}\to \CP^n$. 
We claim that $\wt s$ extends to a \holo\ spray $E\to\CP^n$. Indeed, writing 
$Z=(Z_0:Z')$ with $z=Z'/Z_0$, we see that $s$ has the following expression 
in the homogeneous coordinates $Z\in U_0 \subset\CP^n$:
$$ 
   s(Z,t) = [1 \colon s(Z'/Z_0,t)] = [1\colon Z'/Z_0+ t f(Z/Z_0)] = 
   [Z_0 \colon Z'+  t Z_0 f(Z/Z_0)].  
$$
Hence we get for $k=1,\ldots, n$ and $Z\in U_0\cap U_k$ 
$$
   \wt s \,\theta_k^{-1}(Z,t) = s \,\theta_{0k}(Z,t) 
   = s(Z,t(Z_0/Z_k)^m) = [Z_0 \colon Z'+ t \, Z_0^{m+1}Z_k^{-m} f(Z/Z_0)].
$$
By the choice of $m$ the function $Z_0^{m+1}Z_k^{-m} f(Z/Z_0)$ vanishes
on $\{Z_0=0\}\cap U_k$ and hence $\wt s\,\theta_k^{-1}$ is \holo\ 
on $U_k$. 

This shows that $\wt s \colon E\to \CP^n$ is a spray satisfying 
$\wt s((U_0 \bs A)\times \C) \subset U_0\bs A$ and 
$\wt s(Z,t)=Z$ for all $Z\in A\cup \Lambda$ and $t\in \C$.
For each $Z\in U_0\bs A$ we can find finitely many sprays of this kind
(corresponding to $n$ linearly independent directions in $\C^n$) which 
together dominate at $Z$ and hence at every point in a Zariski open set
containing $Z$. Repeating the construction at other points (and for 
different choices of the hyperplane $\Lambda$) we obtain a finite 
dominating family of algebraic sprays on $\CP_n\bs A$. 

The same proof applies to  complex Grassmanians $Y=G_{k,n}$ since these can 
be covered by finitely many Zariski open \nbd s $U_j\simeq \C^{k(n-k)}$.
\endpr

\demo Proof of Proposition 1.3:
(Compare with Localization Lemma 3.5.B.\ in [G].)
The proof is essentially the same as the proof of Proposition 1.2 (a).
We shall repeatedly use the fact that for every closed algebraic 
subvariety $A\subset Y$ and every point $y\in Y\bs A$ there exists an 
algebraic hypersurface $\Lambda\subset Y$ such that $A\subset \Lambda$ 
but $y\notin\Lambda$. 

Fix a point $y_0\in Y$ and let $U\subset Y$ be a Zariski open \nbd\ of $y_0$ 
with finitely many algebraic sprays $s_j\colon E_j\to Y$ ($j=1,\ldots,k$)
which together dominate at $y_0$. Replacing $U$ by a smaller Zariski open
\nbd\ of $y_0$ we may assume that $\Lambda=Y\bs U$ is an algebraic hypersurface 
in $Y$ and the bundle $E_j|_U\to U$ is algebraically trivial for each $j$. 
Composing an algebraic trivialization of $E_j|_U$ with the spray $s_j$ 
we may therefore assume that $s_j$ is defined on the product bundle 
$U\times\C^{N_j}$ and has values in $Y$. To remove the singularities 
of $s_j$ along $\Lambda$ we replace the product bundle by 
$N_j [\Lambda]^{-m_j}$ for a sufficiently large $m_j\in \N$. 
This gives finitely many sprays on $Y$ which together dominate at $y_0$ 
and hence over a Zariski open \nbd\ of $y_0$. Finitely many such 
collections then dominate on $Y$.
\endpr

\demo Remark: The idea of using the bundles $E=N L^{\otimes m}$ 
where $L=[\Lambda]^{-1}$ can be found in sections 3.5.B.\ and 3.5.C.\ 
of [G]. However, our conclusion is different from the one in [G] where the 
author claimed the existence of a dominating spray on $Y$ obtained by 
composing the individual non-dominating sprays (see the Appendix below). 
The bundles $L$ and $E=NL^{\otimes m}$ used above are strictly negative 
and do not admit any nontrivial \holo\ sections. Hence there exist no nontrivial 
vector bundle maps from any trivial bundle to $E$ (since such a map would take a 
certain constant section to a nontrivial \holo\ section of $E$). Thus we 
are unable to replace the collection of sprays obtained above by a single 
dominating spray using Lemma 2.2.
\endpr

%
%
Proposition 1.5 follows at once from the following more general
result.

\proclaim 5.1 Lemma: Let $A$ be a closed algebraic subvariety
of codimension at least two in a (quasi-) projective algebraic manifold $Y$. 
Suppose that each point $y_0\in Y\bs A$ has a Zariski open \nbd\ $U\subset Y$ 
and an algebraic spray $s\colon E\to Y$, defined on a vector bundle $p\colon E\to U$, 
such that $s$ is dominating at $y_0$ and $s^{-1}(A) \subset E$ contains 
no hypersurfaces. Then $Y\bs A$ is subelliptic.

Indeed, if $s\colon E\to Y$ is a submersive algebraic spray then the codimension 
of $s^{-1}(A)$ in $E$ is the same as the codimension of $A$ in $Y$, and
hence Proposition 1.5 follows.

\demo Proof of Lemma 5.1: 
After removing an algebraic hypersurface which does not contain $y_0$ we may 
assume that $E|_U=U\times \C^N$ and that each fiber 
$\wt A_y=\{t\in\C^N\colon s(y,t)\in A\}$ of $\wt A := s^{-1}(A)$ 
has codimension at least two in $\C^N$. Note that $0\notin \wt A_y$ 
for $y\in U\bs A$. 

Let  $t=(t_1,\ldots,t_N)\in\C^N$. For each $k=0,1,\ldots,N$ we set 
$\C_{(k)}=\C^k\times \{0\}^{N-k}$. 
Let $\pi_k\colon U\times \C_{(k)}\to U\times\C_{(k-1)}$
denote the projection $\pi_k(y,t_1,\ldots,t_k)=(y,t_1,\ldots,t_{k-1})$.
After a linear change of coordinates on $\C^N$ and removing 
another algebraic hypersurface from $U$ we may assume that 

\item{(i)} for each $k=1,\ldots, N$ the set $A_{(k)}=(U\times \C_{(k)}) \cap \wt A$ 
is a subvariety of $U\times \C_{(k)}$ with fibers of codimension at least two
(in particular $A_{(1)}=\emptyset$), and 

\item{(ii)} $\pi_k(A_{(k)}) \subset U\times \C_{(k-1)}$ is an algebraic subvariety 
of $U\times \C_{(k-1)}$ which does not contain the point $(y,0,\ldots,0)$
for any $y\in U$. (Note that $\pi_k(A_{(k)})$ contains $A_{(k-1)}$ 
but it may be larger.) 

\ni Condition (ii) insures that for each $k=2,\ldots,N$ there exists an algebraic 
function $p_k$ on $U\times \C_{(k-1)}$ which vanishes on $\pi_k(A_{(k)})$ and 
satisfies $p_k(y_0,0,\ldots,0)\ne 0$. Consider the map 
$g\colon U\times \C^N\to U\times \C^N$,  
$$ 
	g(y,t)= \bigl( y,t_1,p_2(y,t_1)t_2,\ldots, 
	p_N(y,t_1,\ldots,t_{N-1})t_N \bigr). 
$$
Clearly $g(y,0)=(y,0)$, the map $t\to g(y_0,t)$ is nondegenerate at $t=0$,
and the image of $g$ avoids $\wt A$. Thus $\sigma=s\circ g\colon U\times \C^N\to Y$ 
is an algebraic spray which is dominating at $y_0$ and satisfies 
$\sigma ((U\bs A)\times\C^N) \subset Y\bs A$. Since $y_0$ was an arbitrary point 
of $Y\bs A$, the subellipticity of $Y\bs A$ now follows from Proposition 1.3.
\endpr

%
%
\demo Proof of Proposition 1.6: Let $s\colon E\to Y$ be a dominating 
spray on $Y$ defined on a vector bundle $p\colon E\to Y$.
Denote by $\wt E=\pi^*(E)\to\wt Y$ the pull-back of $E$ by 
the map $\pi \colon \wt Y\to Y$. Explicitly we have
$$  
	\wt E=\{(\wt y,e)\colon \wt y\in \wt Y,\ e\in E,\ \pi(\wt y)=p(e)\}.
$$
Let $\sigma\colon \wt E\to Y$ be defined by $\sigma(\wt y,e)=s(y,e)$
where $y=\pi(\wt y)\in Y$. Fix $\wt y\in \wt Y$. Since the fiber $\wt E_{\wt y}$ 
is simply connected and $\pi$ is a \holo\ covering, the map 
$\sigma(\wt y,\cdotp)\colon \wt E_{\wt y}\to Y$ has a unique holomorphic 
lifting ${\wt s}(\wt y, \cdotp) \colon \wt E_{\wt y} \to \wt Y$ (i.e., 
$\pi( \wt s(\wt y,e)) =\sigma(\wt y,e)$) with ${\wt s}(\wt y, 0)=\wt y$. 
Clearly $\wt s\colon \wt E\to \wt Y$ is a dominating spray on $\wt Y$
and hence $\wt Y$ is elliptic. A similar argument works for  
families of sprays, thereby showing that subellipticity of $Y$ implies 
that of $\wt Y$.

%
%
%
%
\beginsection 6. Removing intersections with complex subvarieties.

Let $X$ and $Y$ be complex manifolds and $A\subset Y$ a closed complex
subvariety. Given a map $f\colon X \to Y$ we write
$f^{-1}(A)=\{x\in X\colon f(x)\in A\}$ and call it the 
{\it intersection set of $f$ with $A$}. If $A$ is a hypersurface
or, more generally, an effective divisor in $Y$, there is a well defined
pull-back divisor $f^*(A)=\sum m_j V_j$ in $X$, where each $V_j$ is an 
irreducible component of the hypersurface $f^{-1}(A)$ and $m_j\in \N$ is 
its multiplicity. The following questions have been 
studied by many authors:

\ni \sl To what extent can the preimage $f^{-1}(A)$ resp.\ $f^*(A)$ be 
prescribed~?  How large is the set of all \holo\ maps $f\colon X\to Y$ 
with the given preimage $f^{-1}(A)$ (resp.\ $f^*(A)$)~?
\rm\medskip

In the simplest case when $X=\C$ and $A$ consists of $d$ points in the
Riemann sphere $Y= \CP^1$ the answer changes when passing from $d=2$ to $d=3$: 
One can prescribe the pull-back of any two points in $\CP^1$ by a \holo\ map 
$f\colon \C\to\CP^1$ (and there are infinitely many such maps), but when 
$d\ge 3$ the pull-back divisor $f^*A$ completely determines the map $f$. 
Similar situation occurs when $A$ consists of $d$ hyperplanes in general 
position in $Y=\CP^n$: we have flexibility up to $d=n+1$ (Corollary 6.2 (b)) 
and rigidity for $d\ge n+2$ (due to hyperbolicity of $\CP^n\bs A$).

The following result shows that subellipticity of the complement $Y\bs A$ 
implies the validity of the Oka principle for maps $f\colon X\to Y$ with the 
given preimage $f^{-1}(A)$ (or pull-back $f^*A$ when $A$ is a divisor). 

\proclaim 6.1 Theorem: Let $A$ be a closed complex subvariety of a complex
manifold $Y$ such that $Y\bs A$ is subelliptic (Definition 2). 
If $X$ is a Stein manifold, $K$ is a compact holomorphically convex 
subset of $X$ and $f\colon X\to Y$ is a continuous map which is \holo\ in an open 
set containing $f^{-1}(A) \cup K$ then for any $r\in\N$ there exist an open set 
$U\supset f^{-1}(A)\cup K$ and a homotopy $f_t\colon X\to Y$
$(t\in[0,1])$ of continuous maps such that $f_0=f$, $f_t$ is \holo\ in $U$ 
and tangent to $f$ to order $r$ along $f_t^{-1}(A)=f^{-1}(A)$ for each 
$t\in [0,1]$, and $f_1$ is \holo\ on $X$.

\medskip\ni\ \bf 6.2 Corollary. \sl  
The conclusion of Theorem 6.1 holds in each of the following cases:
\item{(a)} $Y$ is an affine space $\C^n$, a projective space $\CP^n$
or a complex Grassmanian and $A\subset Y$ is an algebraic subvariety of 
codimension at least two.
\item{(b)} $Y=\CP^n$ and $A$ consists of at most $n+1$ complex hyperplanes 
in general position.
\item{(c)} A complex Lie group acts transitively on $Y\bs A$. 
\rm

In any of these cases $Y\bs A$ is subelliptic by the results in section 1. 
(Note that (b) is a special case of (c).) 

Using Theorem 6.1 we shall also prove the following result which is a version
of the {\it Oka principle for removing of intersections}.

\proclaim 6.3 Theorem: 
Assume that $f\colon X\to Y$ is a \holo\ map, $A$ is a complex subvariety of $Y$
and $f^{-1}(A)=X_0 \cup X_1$, where $X_0$ and $X_1$ are unions of connected components 
of $f^{-1}(A)$ and $X_0\cap X_1=\emptyset$. 
Assume that $X$ is Stein and the manifolds $Y$ and $Y\bs A$ are subelliptic. 
If there exists a homotopy $\wt f_t\colon X\to Y$ $(t\in [0,1])$ 
of continous maps satisfying $\wt f_0=f$, $\wt f_1^{-1}(A)=X_0$, and $\wt f_t|_U=f_t|_U$ 
for some open set $U\supset X_0$ and for all $t\in [0,1]$, then for each 
$r\in\N$ there exists a homotopy of \holo\ maps $f_t\colon X\to Y$ such that 
$f=f_0$, $f_1^{-1}(A)=X_0$, and for each $t\in [0,1]$ the map $f_t$ 
agrees to order $r$ with $f$ along $X_0$ (which is a union of connected 
components of $f_t^{-1}(A)$).

When $A=\{0\}\subset Y=\C^d$, Theorem 6.3 coincides with the main result 
of [FR] on {\it \holo\ complete intersections}. When $Y=\C^d$ and $Y\bs A$ is elliptic 
this is Theorem 1.3 in [F1].  Theorem 6.3 applies if $Y$ is any of the manifolds 
$\C^n$, $\CP^n$ or a complex Grassmanian (these are complex homogeneous and therefore 
elliptic) and $A\subset Y$ is an algebraic subvariety of codimension at least two 
($Y\bs A$ is then subelliptic by Proposition 1.2).

\demo Example: For each $n\ge 1$ there exists a discrete set 
$A\subset\C^n$ for which Theorem 6.3 fails. To see this,
we choose a discrete set $D\subset \C^n$ which is 
{\it unavoidable} in the sense that every entire map 
$F\colon \C^n\to\C^n\bs D$ has rank $<n$ at each point 
(see [RR]). Choose a point $p\in \C^n\bs D$ and set 
$A=D\cup\{p\}$. Take $X=\C^n$, $f=Id \colon\C^n\to\C^n$,
$X_0=\{p\}$ and $X_1=D$. Then the conditions of Theorem 6.3 are 
satisfied but the conclusion fails (since the rank condition for 
holomorphic maps $F\colon \C^n\to \C^n\bs D$ implies that $F^{-1}(p)$ 
contains no isolated points and hence $X_0=\{p\}$ cannot be a
connected component of $F^{-1}(p)$ for any such $F$).
\endpr

Theorem 6.1 is a special case of the following result.

\proclaim 6.4 Theorem: Let $h\colon Z\to X$ be a \holo\ submersion
onto a Stein manifold $X$. Suppose that $Z_0 \subset Z$ is a closed complex 
subvariety of $Z$, $f\colon X\to Z$ is a continuous section and 
$X_0=\{x\in X\colon f(x)\in Z_0\}$. Assume that $f$ is \holo\ in an 
open \nbd\ of $X_0\cup K$ where $K$ is a compact holomorphically convex 
subset of $X$. If the restricted submersion $h\colon Z\bs Z_0 \to X$ is 
subelliptic over $X\bs X_0$ then for each $r\in\N$ there is a homotopy 
$f_t\colon X\to Z$ of continuous sections of $h$ such that $f_0=f$, $f_1$ is 
\holo\ on $X$, and for each $t\in[0,1]$ the section $f_t$ is \holo\ in a 
\nbd\ of $X_0\cup K$, tangent to $f$ to order $r$ along $X_0$ and  
satisfies $\{x\in X\colon f_t(x)\in Z_0\} = X_0$. The analogous result holds
for families of sections.

Indeed we obtain Theorem 6.1 by taking $h \colon Z=X\times Y \to X$ to be
the projection $h(x,y)=x$ and $Z_0=X\times A$. Then 
$h\colon Z\bs Z_0=X\times (Y\bs A)\to X$ is a subelliptic submersion 
when the fiber $Y\bs A$ is subelliptic.

\demo Proof of Theorem 6.4:  We shall follow the proof of Theorem 1.4 in [FP3] 
with some modifications which we shall explain. (The proof is essentially the same 
as the proof of Theorem 1.1 in section 4 above, except for the additional 
interpolation condition on the subvariety $X_0$.) We inductively construct a 
sequence of deformations of the given section which are holomorphic on increasingly 
large open sets in $X$ containing $X_0$ (and exhausting $X$) while at the same 
time paying attention not to introduce any additional intersection points of the
section with the subvariety $Z_0$ (other than $X_0$ where the section is kept 
fixed through the entire process). To insure that no additional intersections 
appear in small \nbd s of $X_0$ we keep the sections tangent to $f$ to a very 
high order along $X_0$ (this will be measured by a suitable coherent sheaf of 
ideals on $X$). Away from $X_0$ we shall perform the modification procedure 
using the restricted submersion $h\colon Z\bs Z_0 \to X$, thereby insuring that 
the sections remain in $Z\bs Z_0$ over $X\bs X_0$.

\proclaim Definition 4:
Let $\cS\subset \cO_X$ be a coherent analytic sheaf of ideals on $X$
and $X_0=\{x\in X\colon \cS_x\ne \cO_{X,x}\}$. 
We say that local \holo\ sections $f_0$ and $f_1$ of $h\colon Z\to X$ 
at a point $x\in X_0$ are $\cS$-tangent at $x$ (denoted $\d_x(f_0,f_1)\in\cS_x$) 
if there exists a local holomorphic chart $\phi$ on $Z$ at $z$ such that
the germ at $x$ of every component of the map $\phi f_0-\phi f_1$ belongs 
to ${\cS}_x$. If $f_0$ and $f_1$ are \holo\ in an open set $U\supset X_0$ 
and $\cS$-tangent at each $x\in X_0$, we say that $f_0$ and $f_1$ are 
$\cS$-tangent and write $\d(f_0,f_1)\in \cS$.

$\cS$-tangency is clearly independent of the choice of local charts on $Z$. 

We now define a sheaf of ideals $\cR\subset \cO_X$ which measures the order of 
contact of a section $f \colon X\to Z$ with a subvariety $Z_0 \subset Z$ along 
$X_0=\{x\in X\colon f(x)\in Z_0\}$. Assume that $f$ is \holo\ in a \nbd\ of $X_0$. 
Fix $x\in X_0$ and let $z=f(x)\in Z_0$. Let $g_1,\ldots,g_k$ be 
\holo\ functions in a \nbd\ $V\subset Z$ of $z$ which generate the sheaf of 
ideals of $Z_0 \cap V$ at each point of $V$. Let $\cR_V$ denote the sheaf 
of ideals in $\cO_X|_V$ generated by the functions $g_j \circ f$, $1\le j\le k$. 
(For $x\in V\bs X_0$ we have $\cR_x=\cO_x$.) It is easily seen that $\cR_V$ 
does not depend on the choice of the local generators $g_j$ and hence we 
obtain a coherent analytic sheaf of ideals $\cR \subset \cO_X$ with support $X_0$. 

Fix an integer $r\in \N$ and let $\cS=\cR\cdotp \cJ^r$, where $\cJ \subset \cO_X$ 
is the sheaf of ideals of $X_0$. Let $d$ be a metric on $Z$.
The following lemma was proved in [Fo1, sect.\ 3].

\proclaim 6.5 Lemma: 
Let $f\colon X\to Z$ be a continuous section of $h\colon Z\to X$ which 
is \holo\ in an open set $U$ containing $X_0=\{x\in X\colon f(x)\in Z_0\}$. 
If $g\colon U\to Z$ is a \holo\ section of $h\colon Z|_U \to U$ 
satisfying $\d(f,g)\in \cS$ then there is an open set $V \supset X_0$ 
such that $\{x\in V\colon g(x)\in Z_0\} = X_0$. Furthermore if 
$K\ss K'$ are compact sets in $U$ then there is an $\e>0$ 
such that for any $g$ as above satisfying $d(f(x),g(x))<\e$
$(x\in K')$ we have $\{x\in K\colon g(x)\in Z_0\}=X_0\cap K$.

Everything is now ready for us to explain the modification problem
which is the main ingredient in the proof of Theorem 6.4. 

\ni\bf The modification problem. \rm (Assumptions as in Theorem 6.4.)
Let $B\subset X\bs X_0$ be a compact set such that $(K,B)$ is a Cartan pair 
in $X$ (section 4) and $K\cup B$ is holomorphically convex in $X$.
Let $A=X_0\cup K$ and let $a\colon \wt A\to Z$, $b\colon \wt B\to Z\bs Z_0$ 
be \holo\ sections of $h$ in open \nbd s $\wt A\supset A$ resp.\ $\wt B\supset B$ 
such that $\{x\in\wt A\colon a(x)\in Z_0\}=X_0$. Suppose furthermore that 
$b_t \colon \wt C\to Z\bs Z_0$ ($t\in[0,1]$) is a family of \holo\ sections 
over $\wt C=\wt A\cap \wt B \subset X\bs X_0$ such that $b_0=b|_{\wt C}$ and 
$b_1=a|_{\wt C}$. From these data we must construct a homotopy $(\wt a_t,\wt b_t)$ 
$(t\in[0,1])$ of \holo\ sections in smaller open \nbd s $A'\supset A$ resp.\ 
$B'\supset B$ satisfying the following properties: 

\smallskip
\item{(i)}  $\wt a_0=a|_{A'}$, $\wt b_0=b|_{B'}$, 
\item{(ii)} $\wt a_1=\wt b_1$ on $A'\cap B'$, and hence this pair defines  
a \holo\ section $\wt a$ on $A'\cup B'$, 
\item{(iii)} for each $t\in [0,1]$ we have $\d(\wt a_t,a)\in\cS$, 
$\{x\in A'\colon \wt a_t(x)\in Z_0\}=X_0$, and $\wt a_t|_K$ is uniformly close 
to $a|_K$, and
\item{(iv)} $\wt b_t(x)\in Z\bs Z_0$ for all $x\in B'$ and $t\in [0,1]$.

\smallskip
The final section $\wt a$ is \holo\ on $A'\cup B' \supset X_0\cup K\cup B$, 
it is $\cS$-tangent to $a$ along $X_0$ and its graph intersects $Z_0$ precisely 
over $X_0$ as required. This will complete the induction step. 
To prove the parametric version of Theorems 6.4 and 6.1 one must also consider 
the analogous modification problem for continuous families of sections 
$\{(a_p,b_p)\colon p\in P\}$ on a \nbd\ of $A$ resp.\ $B$, where 
$P$ is a compact Hausdorff space. This extension presents no additional 
difficulties and we refer to [FP1] for the details.

To solve the above modification problem we proceed as in section 4 above.
By Theorem 3.1 (the Oka-Weil theorem) we can deform $b$ through a homotopy of 
\holo\ sections of $Z\bs Z_0$ over a \nbd\ of $B$ to another section 
(still denoted $b$) which approximates $a$ uniformly on a \nbd\ of $C=A\cap B$.
The remaining problem is to patch $a$ and $b$. This will be done as in section 4,
but with a couple of modifications which we now explain. As in sect.\ 4
we find open \nbd s $A' \supset A$, $B'\supset B$ of $A$ resp.\ $B$  and 
\item{(a)} a local $h$-spray $s_1 \colon A'\times \B^n(\eta) \to Z$ with
$s_1(x,0)=a(x)$ for $x\in A'$,
\item{(b)} a global $h$-spray $s_2\colon B'\times\C^n\to Z\bs Z_0$ with
$s_2(x,0)=b(x)$ for $x\in B'$, and 
\item{(c)} a transition map $\psi\colon C'\times \B^n(\eta)\to \C^n$ satisfying 
$$ 
	s_2(x,\psi(x,t))=s_1(x,t), \quad (x\in C'=A'\cap B',\ t\in\B^n(\eta)).
$$

The only addition is that we build $\cS$-tangency into the construction of the
local spray $s_1$ as follows. First we choose a preliminary local dominating 
$h$-spray $s'_1\colon A'\times\B^k(\eta') \to Z$ with $s'_1(x,0)=a(x)$. 
By Cartan's Theorem A [GR] there exist finitely many 
global sections $h_1,\ldots, h_m$ of the sheaf 
$\cS$ on $X$ such that $X_0=\{x\in X\colon h_j(x)=0,\ 1\le j\le m\}$.  
(The $h_j$'s need not generate $\cS$.) Define $\tau\colon A'\times (\C^k)^m\to \C^k$ 
by $\tau(x,t_1,\ldots,t_m)=\sum_{j=1}^m h_j(x)t_j$, where $t_j\in\C^k$
for each $j$. Set $n=mk$ and $t=(t_1,\ldots,t_n)\in\C^n$.
Then for suitably small $\eta>0$ and $A'\supset A$ the map
$$ 
	s_1(x,t)=s'_1(x,\tau(x,t)), \quad (x\in A',\ t\in  \B^n(\eta)) 
$$
is a local $h$-spray which is dominating on $A'\bs X_0$, and any section of the 
form $a_{\a}(x)=s_1(x,\a(x))$ (where $\a\colon A' \to \B^n(\eta)$ is a \holo\ map) 
is $\cS$-tangent to $a$ along $X_0$. If $\eta>0$ and $A'\supset A$ are chosen suficiently 
small then Lemma 6.5 insures that the graph of $a_\a$ intersects $Z_0$ precisely over $X_0$.

After shrinking the sets $A'\supset A$ and $B'\supset B$ we obtain by Proposition 4.1 
in [FP3] a pair of \holo\ maps $\a\colon A'\to \B^n(\eta)$,
$\b\colon B'\to\C^n$ such that $\psi(x,\a(x))=\b(x)$ for $x\in C'=A'\cap B'$.
The \holo\ homotopies
$$ 
   \wt a_t(x)=s_1(x,t\a(x))\ \ (x\in A'),\quad 
   \wt b_t(x)=s_2(x,t\b(x))\ \ (x\in B') 
$$ 
for $t\in [0,1]$ then solve the modification problem. In fact for $t=1$ 
and $x\in C'$ we have 
$$ \wt b_1(x)=s_2(x,\b(x))=s_2(x,\psi(x,\a(x)))=s_1(x,\a(x))=\wt a_1(x), $$ 
and the other requirements are easily verified.

Using the solution of this modification problem we prove Theorem 6.4 by 
following the globalization scheme in the proof of Theorem 1.4 in [FP3]. 
More precisely, we follow the second approach in [FP3, pp.\ 65-66] whose main 
advantage is that the patching of (families of) sections is performed only on 
sets $C= A\cap B\subset X\bs X_0$ and hence no special condition 
on the submersion $h$ is required over $X_0$. (In [F2] it is shown that 
we may even allow $h$ to have ramification points, provided that these 
project by $h$ into the subvariety $X_0$.) With the same tools one can obtain
the extension of Theorem 6.4 to continuous families of sections 
$\{f_p\colon p\in P\}$ with the parameter in a compact Hausdorff space $P$ 
(see [FP2, FP3]).
\endpr

\demo Proof of Theorem 6.3: We replace maps $X\to Y$ by 
sections of $Z=X\times Y\to X$ without changing the notation. 
By hypothesis there is a continuous homotopy $\wt f_t \colon X\to Z$ 
($t\in [0,1]$), with $\wt f_0=f$, which is fixed near $X_0$
and satisfies $\{x\in X\colon \wt f_1(x)\in Z_0\}=X_0$. We now apply
Theorem 6.4, with $\wt f_1$ as the initial section, to obtain a homotopy
$\wt f_t \colon X\to Z$ ($t\in [1,2]$), where the final section 
$\wt f_2$ is \holo\ on $X$ and satisfies 
$\{x\in X\colon \wt f_2(x)\in Z_0\}=X_0$. We rescale the parameter 
interval $[0,2]$ back to $[0,1]$. Since $Z\to X$ is assumed to be subelliptic 
over $X\bs X_0$, we can apply Theorem 6.4 
to the homotopy $\{\wt f_t \colon t\in [0,1]\}$, 
with the sheaf $\cS$ defined above (for some fixed $r\in\N$), to obtain 
a two-parameter homotopy $h_{t,s}\colon X\to Z$ ($t,s\in [0,1]$) of 
continuous sections which are \holo\ in a \nbd\ $V\supset X_0$ 
(independent of $t,s$) and satisfy:

\smallskip
\item{(i)} $h_{t,0}=\wt f_t$ for all $t\in [0,1]$,
\item{(ii)} $h_{0,s}=\wt f_0=f_0$ and $h_{1,s}=\wt f_1$ for all
$s\in [0,1]$,
\item{(iii)} $\d(h_{t,0},h_{t,s}) \in \cS$ for all
$s,t\in [0,1]$, and
\item{(iv)} the map $f_t:=h_{t,1}$ is \holo\ on $X$
for each $t\in [0,1]$.

\smallskip 
It follows from (iii) and Lemma 6.5 that there is a \nbd\ $U\subset V$ 
of $X_0$ such that $\{x\in U\colon h_{t,s}(x)\in Z_0\} =X_0$ for all 
$s,t\in [0,1]$. The homotopy $\{f_t\colon t\in [0,1]\}$ defined by
(iv) above then satisfies the conclusion of Theorem 6.3. 
\endpr

\medskip\ni \bf Appendix: Remarks on the paper [G]. \rm

We wish to point out certain inconsistencies in 
section 3 of the paper [G], in particular those 
concerning the notion $\Ell_\infty$. 

\ni 1.\ In the holomorphic category the property $\Ell_\infty$ for a complex 
manifold $Y$ means the validity of a certain strong form of the Oka principle 
for maps from all Stein manifolds to $Y$ (section 3.1.\ in [G]). However, 
in [G, sect.\ 3.5] an algebraic manifold $Y$ is said to be 
{\it algebraically $Ell_\infty$} if it admits a dominating algebraic spray 
$s\colon E\to Y$, defined on an algebraic vector bundle $E\to Y$. 
The conclusion of Localization Lemma 3.5.B. in [G] is that a manifold $Y$
satisfying the hypothesis of that lemma is algebraically $\Ell_\infty$.
However, the proof offered there only gives finitely many algebraic sprays 
$s_i\colon E_i\to Y$ which together dominate at each point $y\in Y$ 
(thus showing that $Y$ is subelliptic in our sense) and concludes with the sentence: 
{\it Then a composition of finitely many such $s_i$...gives us the 
desired dominating spray over $Y$.} Since the bundles $E_i$ are not 
necessarily trivial (in the proof of Proposition 1.2 they are in fact
negative and do not admit any nontrivial holomorphic sections), the sprays 
$s_i$ cannot be pulled back to nondegenerate sprays on trivial bundles 
over $Y$ which would be necessary in order to use Lemma 2.4. 

\ni 2. A similar remark applies to section 3.5.C.\ in [G] which claims 
the implication (4.5) for any Zariski closed subset $A$ of codimension at 
least two in an algebraic manifold $Y$. Being unable to prove 
this we offer Proposition 1.5 (and Lemma 5.1) above as a  replacement.

\ni 3. Another inconsistency can be found in [G, sect.\ 3.2.A'].
The question considered in [G, sect.\ 3.2] is whether the $\Ell_\infty$ 
property of a complex manifold $Y$ implies the existence of a dominating
spray on $Y$. In [G, sec.\ 3.2.A'] it is claimed that this is the case 
if $Y$ is a projective variety which admits a `sufficiently negative' 
vector bundle $E\to Y$ whose rank $N$ is large compared to $\dim Y$. 
The idea in [G] is to first construct a `local spray' $s_0$ on $Y$, 
defined in a small tubular \nbd\ of the zero section $Y_0 \subset E$ in $E$, 
and subsequently Runge approximate $s_0$ by a global \holo\ map 
$s\colon E\to Y$ (which is then a dominating spray on $Y$).  
Here the author refers to the assumed axiom $\Ell_\infty$ 
which pertains to maps from Stein manifolds to $Y$. The problem 
is that the manifold $E$ is not Stein (it contains the 
compact complex submanifold $Y_0$). When the bundle $E\to Y$
is negative, $E$ is holomorphically convex and admits an exhaustion 
function which is zero on $Y_0$ and \spsh\ on $E\bs Y_0$. 
In order to make the above conclusion valid one would have to 
change the axiom $\Ell_\infty$ so that it would pertain to maps from 
all holomorphically convex manifolds (and not only Stein manifolds)
to the given manifold $Y$. Although it seems likely that subellipticity 
of $Y$ implies a suitable version of the Oka principle for maps $X\to Y$ for 
any holomorphically convex manifold $X$, no such result has been 
proved yet.

\ni 4. Section 3 in [G] contains several further results which we have not 
been able to fully understand and justify. For instance, in sect.\ 3.5.C.\ 
of [G] one finds the following statement: {\it If a Zariski closed subset 
$Y_0\subset Y$ in an algebraic manifold $Y$ satisfies ${\rm codim} Y_0\ge 2$ then} 
$$
   (\Ell_\infty\ {\rm for}\ Y) \Rightarrow (\Ell_\infty\ {\rm for}\ Y'=Y\bs Y_0).					 
$$
Here $\Ell_\infty$ property for $Y$ means the existence of a dominating algebraic 
spray on $Y$ according to the first sentence in [G, 3.5.A]. Although this may 
well be the case, we are unable to justify the argument which is supposed to bring 
a spray on $Y$ in `general position' with respect to the subvariety $Y_0$ 
(see the discussion in sect.\ 3.5.C' in [G], and compare with our Proposition 1.5
and Lemma 5.1.)

%
%
\medskip \ni\bf Acknowledgements. \rm
This article is based to a large extent on the joint papers 
[FP1, FP2, FP3] of the author with J.\ Prezelj. I wish to thank her for 
the continuing interest and stimulating discussions on this subject. 
I thank J.\ Winkelmann for several useful observations regarding 
the existence of sprays. This research was supported in part by a grant 
from the Ministry of Science of the Republic of Slovenia and by an NSF grant.

%
%
%
%
\medskip\ni\bf References. \rm

\ii{[C]} H.\ Cartan: Espaces fibr\'es analytiques.
Symposium Internat.\ de to\-pologia algebraica, Mexico, 97--121 (1958).
(Also in Oeuvres, vol.\ 2, Springer, New York, 1979.)

\ii{[GH]} P.\ Griffiths, J.\ Harris:
{\it Principles of Algebraic Geometry}.
John Wiley \& Sons, New York, 1978.

\ii{[FR]} O.\ Forster and K.\ J.\ Ramspott:
Analytische Modulgarben und Endromisb\"undel.
Invent.\ Math.\ {\bf 2}, 145--170 (1966).

\ii{[F1]} F.\ Forstneri\v c: On complete intersections. 
Ann. Inst. Fourier {\bf 51} (2001), 497--512.

\ii{[F2]} F.\ Forstneri\v c:
The Oka principle for multivalued sections of ramified mappings. 
Preprint, 2001.

\ii{[FP1]} F.\ Forstneri\v c and J.\ Prezelj:
Oka's principle for holomorphic fiber bundles with sprays.
Math.\ Ann.\ {\bf 317} (2000), 117-154.

\ii{[FP2]} F.\ Forstneri\v c and J.\ Prezelj:
Oka's principle for holomorphic submersions with sprays.
Math.\ Ann., to appear.

\ii{[FP3]} F.\ Forstneri\v c and J.\ Prezelj: 
Extending holomorphic sections from complex subvarieties.
Math.\ Z.\ {\bf 236} (2001), 43--68.

\ii{[Gr1]} H.\ Grauert:
Holomorphe Funktionen mit Werten in komplexen Lieschen Gruppen.
Math.\ Ann.\ {\bf 133}, 450--472 (1957).

\ii{[Gr2]} H.\ Grauert: Analytische Faserungen \"uber
holomorph-vollst\"andi\-gen R\"aumen.
Ma\-th.\ Ann.\ {\bf 135}, 263--273 (1958).

\ii{[G]} M.\ Gromov:
Oka's principle for holomorphic sections of elliptic bundles.
J.\ Amer.\ Math.\ Soc.\ {\bf 2}, 851-897 (1989).

\ii{[GR]} C.\ Gunning, H.\ Rossi:
Analytic functions of several complex variables.
Prentice--Hall, Englewood Cliffs, 1965.

\ii{[HL]} G.\ Henkin, J.\ Leiterer:
The Oka-Grauert principle without induction over the basis dimension.
Math.\ Ann.\ {\bf 311}, 71--93 (1998).

\ii{[L]} F.\ L\'arusson:
Excision for simplicial sheaves on the Stein site
and Gromov's Oka principle. Preprint, December 2000.

\ii{[O]} K.\ Oka: Sur les fonctions des plusieurs variables. III:
Deuxi\`eme probl\`eme de Cousin.
J.\ Sc.\ Hiroshima Univ.\ {\bf 9} (1939), 7--19.

\ii{[RR]} J.-P.\ Rosay, W.\ Rudin:
Holomorphic maps from $\C^n$ to $\C^n$.
Trans.\ Amer.\ Math.\ Soc.\ {\bf 310}, 47--86 (1988).

\bigskip\medskip
\ni\bf Address: \rm Institute of Mathematics, Physics and Mechanics,
University of Ljubljana, Jadranska 19, 1000 Ljubljana, Slovenia

\bye